\input amstex
\documentstyle{amsppt}
\magnification=\magstep1

\vsize 21.8 true cm \hsize 15.2 true cm \voffset 1 true cm
\lineskip=3pt plus 1pt minus 1pt

\leftheadtext{C. Baak, H. Chu, M. Moslehian} \rightheadtext{$n$-inner product preserving mappings}

\topmatter

\title On the Cauchy--Rassias Inequality and Linear $n$-Inner Product Preserving 
Mappings
\endtitle

\author Choonkil Baak$^1$, Hahng-Yun Chu$^2$ and Mohammad Sal Moslehian$^3$
\endauthor

\address
$^{1, 2}$Department of Mathematics, Chungnam National University, Daejeon
305--764, South Korea; $^3$Department of Mathematics, Ferdowsi University, P.O. Box 1159,
Mashhad 91775, Iran
\endaddress

\email $^1$cgpark\@cnu.ac.kr; $^2$hychu\@math.cnu.ac.kr; $^3$msalm\@math.um.ac.ir
\endemail

\abstract
We prove the Cauchy--Rassias stability of linear $n$-inner product preserving mappings in
$n$-inner product Banach spaces. We apply the Cauchy-Rassias inequality that plays an influencial role in the subject of functional equations. The inequality was introduced for the first time by Th.M.Rassias in his paper entitled: On the stability of the linear mapping in Banach spaces, Proc. Amer. Math.Soc. 72(1978), 297--300.

\endabstract

\thanks Supported by Korea Research Foundation Grant
KRF-2004-041-C00023.
\endthanks

\keywords Cauchy--Rassias stability, linear $n$-inner product preserving mapping,
$n$-inner product Banach space, Hilbert space
\endkeywords

\subjclass Primary 39B52, 46B04, 46B20, 51Kxx
\endsubjclass
\endtopmatter

\document

\baselineskip=0.6 true cm

\head 1.  Introduction \endhead

Let $X$ and  $Y$  be Banach spaces with norms $|| \cdot ||$ and
 $\| \cdot \|$, respectively. Consider $f : X
\rightarrow Y$ to be a mapping such that $f(tx)$ is continuous in
$t \in \Bbb R$ for each fixed $x \in X$. Th.M. Rassias \cite{23}
introduced the following inequality, that we call {\it
Cauchy--Rassias inequality} :  Assume that there exist constants
$\theta \ge 0$ and $p\in [0, 1)$ such that
$$\| f(x+y)-f(x)-f(y)\|\le \theta (||x||^p+||y||^p)\tag *$$  for
all $x, y\in X$. Rassias \cite{23} showed that there exists a
unique $\Bbb R$-linear mapping $T : X \rightarrow Y$ such that
$$\|f(x) -T(x)\|\le \frac{2\theta}{2-2^p}||x||^p$$ for all $x
\in X$. The inequality {\rm ($*$)} has provided a lot of influence
in the development of what we now call {\it Hyers--Ulam--Rassias
stability} of functional equations. Beginning around the year 1980
the topic stability of linear functional equations has been studied by a number of
mathematicians  (see \cite{7}, \cite{8} and \cite{11}--\cite{27}).

Jun and Lee \cite{9} proved the following: Denote by
 $\varphi : X \setminus \{0\} \times X \setminus \{0\}
\rightarrow [0, \infty)$ a function such that
$$ \widetilde{\varphi}(x, y)=\sum_{j=0}^{\infty} \frac{1}{
3^j} \varphi(3^jx, 3^jy)  <  \infty $$ for all $x, y \in X
\setminus \{0\}$. Suppose that $f : X \rightarrow Y$ is a mapping
satisfying $$\|2f(\frac{x+y}{2})-f(x)-f(y)\| \le  \varphi(x, y)$$
for all $x, y \in X \setminus \{0\}$. Then there exists a unique
additive mapping $T : X \rightarrow Y$ such that
$$\|f(x)-f(0)-T(x)\| \le  \frac{1}{3}\big(\widetilde{\varphi}(x, -x) +
\widetilde{\varphi}(-x, 3x)\big) $$ for all $x\in X \setminus \{0\}$.

Let $X$ and $Y$ be complex Hilbert spaces. A mapping $f : X \rightarrow Y$ is called
an {\it inner product preserving mapping } if $f$ satisfies the {\it orthogonality equation }
$$\langle f(x) , f(y) \rangle = \langle x , y \rangle $$
for all $x, y\in X$.
The inner product preserving mapping problem has been investigated in several
papers (see \cite{2}--\cite{4}).

In this paper, we prove the Cauchy--Rassias stability of inner product
preserving mappings in Hilber spaces, introduce the concept of $n$-inner
product Banach space and establish the Cauchy--Rassias stability of linear
 $n$-inner product preserving mappings in $n$-inner product Banach spaces.

\head 2. Cauchy--Rassias stability of inner product preserving mappings in Hilbert spaces
\endhead

In this section, assume that $X$ is a complex Hilbert space with norm $||~ \cdot~ ||$,  and that $Y$ is
a complex Hilbert space with norm  $\|~ \cdot ~\|$.

\proclaim{Theorem 2.1}
 Let $f : X \rightarrow Y$ be a mapping with $f(0) = 0$
 for which there exists a function $\varphi
: X \times X \rightarrow [0, \infty)$  such that
$$\align
\widetilde{\varphi}(x, y) : =  \sum_{j=0}^{\infty} \frac{1}{2^{j}}
\varphi(2^j x, 2^j y) &  < \infty , \tag 2.i \\
 \|f(\mu x + \mu y) -  \mu f(x) - \mu f(y) \| & \le \varphi(x, y) ,  \tag 2.ii \\
| ~ \langle f(x), f(y)\rangle - \langle x, y \rangle ~ | & \le \varphi(x, y)  \tag 2.iii
\endalign$$ for  all $\mu \in \Bbb T^1 : = \{\lambda \in \Bbb C ~ \mid ~ |\lambda | = 1\}$ and all $x, y \in X$.
 Then there exists a unique inner product preserving mapping $U : X \rightarrow Y$ such that
$$\|f(x)- U(x)\| \le \frac{1}{2} \widetilde{\varphi}(x, x)
\tag 2.iv $$
for all $x \in X$.
\endproclaim

\demo{Proof}
By the same reasoning as in the proof of \cite{21, Theorem 2.1}, it follows from {\rm (2.i)} and {\rm (2.ii)} that
there exists a unique $\Bbb C$-linear mapping $U : X \rightarrow
Y$ satisfying {\rm (2.iv)}. The $\Bbb C$-linear mapping $U : X \rightarrow Y$ is given by
$$U(x) = \lim_{l \to \infty} \frac{1}{2^l} f(2^l x) \tag 2.1
$$ for all $x \in X$.

It follows from {\rm (2.iii)} that $$\align  | ~ \langle \frac{1}{2^l}
f(2^l x), \frac{1}{2^l} f(2^l y)\rangle - \langle x, y \rangle ~ | & =
 \frac{1}{2^{2l}} | ~ \langle f(2^l x), f(2^l y)\rangle - \langle 2^l x, 2^l y\rangle ~ | \\ & \le \frac{1}{2^{2l}}
 \varphi(2^l x, 2^l y) \le \frac{1}{2^{l}}
 \varphi(2^l x, 2^l y) ,
\endalign $$
which tends to zero as $l \to \infty$ for all $x, y \in X$  by {\rm (2.i)}. By {\rm (2.1)},
$$ \langle U(x) , U(y) \rangle  = \lim_{l \to \infty} \langle \frac{1}{2^l}
f(2^l x), \frac{1}{2^l}
f(2^l y) \rangle = \langle x, y \rangle $$ for all $x, y \in X$, as desired.
\qed
\enddemo

\proclaim{Corollary 2.2}  Let $f : X \rightarrow Y$ be a
mapping with $f(0) = 0$ for which there exist constants
$\theta \ge 0$ and $p \in [0, 1)$ such that
$$\align \| f(\mu x + \mu y)- \mu f(x) - \mu f(y)\| & \le  \theta
(||x||^p + ||y||^p) , \\
|~ \langle f(x) , f(y) \rangle - \langle x, y \rangle  ~| & \le \theta (||x||^p + ||y||^p)
\endalign$$ for all $\mu \in \Bbb T^1$ and all  $x, y \in X$.
  Then there exists a unique inner product preserving mapping $U :  X
\rightarrow Y$ such that $$\|f(x)- U(x)\| \le \frac{2\theta}{2-2^p} ||x||^p $$
for all $x \in X$.
\endproclaim

\demo{Proof} Define $\varphi(x, y) = \theta (||x||^p + ||y||^p)$
to be Th.M. Rassias upper bound in the Cauchy--Rassias inequality,
and apply Theorem 2.1. \qed
\enddemo

\proclaim{Theorem 2.3}
 Let $f : X \rightarrow Y$ be a
mapping with $f(0)= 0$ for which there exists a function $\varphi
: X \setminus \{0\} \times X \setminus \{0\} \rightarrow [0, \infty)$  such that
$$\align
\widetilde{\varphi}(x, y) : = \sum_{j=0}^{\infty} \frac{1}{3^j}
\varphi(3^j x, 3^j y) & < \infty , \tag 2.v \\
 \|2f(\frac{\mu x + \mu y}{2}) -  \mu  f(x) - \mu f(y)\|
& \le \varphi(x, y) ,  \tag 2.vi \\
| ~ \langle f(x), f(y)\rangle - \langle x, y \rangle ~ | & \le \varphi(x, y) \tag 2.vii
\endalign$$ for all $\mu \in \Bbb T^1$ and all  $x, y \in X\setminus \{ 0 \}$.
 Then there exists a unique inner product preserving mapping $U :
X \rightarrow Y$ such that
$$\|f(x)- U(x)\| \le  \frac{1}{3}\big(\widetilde{\varphi}(x, -x) +
\widetilde{\varphi}(-x, 3x)\big) \tag 2.viii $$
for all $x \in X \setminus \{0\}$.
\endproclaim

\demo{Proof} By the same reasoning as in the proof of \cite{11, Theorem 2.5},
it follows from {\rm (2.v)} and {\rm (2.vi)} that
there exists a unique $\Bbb C$-linear mapping $U : X \rightarrow
Y$ satisfying {\rm (2.viii)}. The $\Bbb C$-linear mapping $U : X
\rightarrow Y$ is given  by
$$U(x) = \lim_{l\to\infty}\frac{1}{3^l} f(3^l x) \tag 2.2 $$
for all $x \in X$.

It follows from {\rm (2.vii)} that $$\align  | ~ \langle \frac{1}{3^l}
f(3^l x), \frac{1}{3^l} f(3^l y)\rangle - \langle x, y \rangle ~ | & =
 \frac{1}{3^{2l}} | ~ \langle f(3^l x), f(3^l y)\rangle - \langle 3^l x, 3^l y\rangle ~ | \\ & \le \frac{1}{3^{2l}}
 \varphi(3^l x, 3^l y) \le \frac{1}{3^{l}}
 \varphi(3^l x, 3^l y) ,
\endalign $$
which tends to zero as $l \to \infty$ for all $x, y \in X$  by {\rm (2.v)}. By {\rm (2.2)},
$$ \langle U(x) , U(y) \rangle  = \lim_{l \to \infty} \langle \frac{1}{3^l}
f(3^l x), \frac{1}{3^l}
f(3^l y) \rangle = \langle x, y \rangle $$ for all $x, y \in X$, as desired.
\qed
\enddemo

\proclaim{Corollary 2.4}  Let $f : X \rightarrow Y$ be a
mapping with $f(0) = 0$ for which there exist constants
$\theta \ge 0$ and $p \in [0, 1)$ such that
$$\align \| 2f(\frac{\mu x + \mu y}{2})- \mu f(x) - \mu f(y)\| & \le  \theta
(||x||^p + ||y||^p) , \\
| ~ \langle f(x), f(y)\rangle - \langle x, y \rangle ~ | & \le \theta (||x||^p + ||y||^p)
\endalign$$ for all $\mu \in \Bbb T^1$ and all  $x, y \in X \setminus \{ 0 \}$.
  Then there exists a unique inner product preserving mapping $U :  X
\rightarrow Y$ such that $$\|f(x)- U(x)\| \le \frac{3+3^p}{3-3^p} \theta ||x||^p $$
for all $x \in X\setminus \{ 0 \}$.
\endproclaim

\demo{Proof} Define $\varphi(x, y) = \theta (||x||^p + ||y||^p)$,
and apply Theorem 2.3. \qed
\enddemo

\proclaim{Theorem 2.5}
 Let $f : X \rightarrow Y$ be a
mapping with $f(0)= 0$ for which there exists a function $\varphi
: X \setminus \{0\} \times X \setminus \{0\} \rightarrow [0, \infty)$  such that
$$\align
 \sum_{j=0}^{\infty} 3^{2j}
\varphi(\frac{x}{3^j}, \frac{y}{3^j}) & < \infty , \tag 2.ix \\
 \|2f(\frac{\mu x + \mu y}{2}) - \mu  f(x) - \mu f(y)\|
& \le \varphi(x, y) ,  \tag 2.x \\
| ~ \langle f(x), f(y)\rangle - \langle x, y \rangle ~ | & \le \varphi(x, y)  \tag 2.xi
\endalign$$ for all $\mu \in \Bbb T^1$ and all  $x, y \in X\setminus \{ 0 \}$.
 Then there exists a unique inner product preserving mapping  $U :
X \rightarrow Y$ such that
$$\|f(x)- U(x)\| \le  \widetilde{\varphi}(\frac{x}{3}, - \frac{x}{3}) +
\widetilde{\varphi}(-\frac{x}{3}, x) \tag 2.xii $$ for all $x \in
X \setminus \{0\}$, where $$ \widetilde{\varphi}(x, y) : =
\sum_{j=0}^{\infty} 3^j \varphi(\frac{x}{3^j}, \frac{y}{3^j})$$
for all $x, y \in X\setminus \{ 0 \}$.
\endproclaim

\demo{Proof} Note that $$3^{j} \varphi(\frac{x}{3^j}, \frac{y}{3^j}) \le
3^{2j} \varphi(\frac{x}{3^j}, \frac{y}{3^j}) \tag 2.3$$
for all $x, y \in X$ and all positive integers $j$.
By the Jun and Lee's theorem  \cite{9, Theorem 6}, it follows from {\rm (2.ix)}, {\rm (2.3)}  and {\rm (2.x)} that
there exists a unique additive mapping $U : X \rightarrow
Y$ satisfying {\rm (2.xii)}. The additive mapping $U : X
\rightarrow Y$ is given  by
$$U(x) = \lim_{l\to\infty} 3^l f(\frac{x}{3^l}) \tag 2.4 $$
for all $x \in X$.
By the same method as in the proof of \cite{11, Theorem 2.5},  one can show that
the additive mapping  $U : X \rightarrow
Y$ is $\Bbb C$-linear.

It follows from {\rm (2.xi)} that $$ | ~ \langle 3^l
f(\frac{x}{3^l}), 3^l f(\frac{y}{3^l})\rangle - \langle x, y \rangle ~ |  =
 3^{2l} | ~ \langle f(\frac{x}{3^l}), f(\frac{y}{3^l})\rangle - \langle \frac{x}{3^l}, \frac{y}{3^l}
\rangle ~ |  \le 3^{2l} \varphi(\frac{x}{3^l}, \frac{y}{3^l}) ,
$$
which tends to zero as $l \to \infty$ for all $x, y \in X$  by {\rm (2.ix)}. By {\rm (2.4)},
$$ \langle U(x) , U(y) \rangle  = \lim_{l \to \infty} \langle 3^l
f(\frac{x}{3^l}), 3^l f(\frac{y}{3^l}) \rangle = \langle x, y \rangle $$ for all $x, y \in X$, as desired.
\qed
\enddemo

\proclaim{Corollary 2.6}  Let $f : X \rightarrow Y$ be a
mapping with $f(0) = 0$ for which there exist constants
$\theta \ge 0$ and $p \in (2, \infty)$ such that
$$\align \| 2f(\frac{\mu x + \mu y}{2})- \mu f(x) - \mu f(y)\| & \le  \theta
(||x||^p + ||y||^p) , \\
| ~ \langle f(x), f(y)\rangle - \langle x, y \rangle ~ | & \le \theta (||x||^p + ||y||^p)
\endalign$$ for  all $\mu \in \Bbb T^1$ and all  $x, y \in X \setminus \{ 0 \}$.
  Then there exists a unique inner product preserving mapping  $U :  X
\rightarrow Y$ such that $$\|f(x)- U(x)\| \le \frac{3^p+3}{3^p-3} \theta ||x||^p $$
for all $x \in X\setminus \{ 0 \}$.
\endproclaim

\demo{Proof} Define $\varphi(x, y) = \theta (||x||^p + ||y||^p)$,
and apply Theorem 2.5. \qed
\enddemo

\head 3. Cauchy--Rassias stability of linear $n$-inner product preserving mappings in $n$-inner product Banach spaces
\endhead

We first recall the notion of $n$-inner product space.

\definition{Definition 3.1 \cite{5}} Let $X$ be a complex linear
space with $\dim X \ge n \ge 2$ and $\langle \cdot, \cdot \mid \cdot,
\cdots, \cdot\rangle : X^{n+1} \to \Bbb C$ be a function. Then $(X,
\langle \cdot, \cdot \mid \cdot, \cdots, \cdot\rangle )$ is called
an {\it $n$-inner product space} if
$$\align
&(\text{\rm nI}_1) \,\ \langle x, x \mid x_2, \cdots, x_n \rangle \ge 0 , \\
&(\text{\rm nI}_2) \,\ \langle x, x \mid x_2, \cdots, x_n \rangle = 0  \text {if and only if} x, x_2, \cdots, x_n
\text{ are linearly
dependent} , \\
&(\text{\rm nI}_3) \,\ \langle x, y \mid x_2, \cdots, x_n \rangle = \overline{\langle y,
 x \mid x_2, \cdots, x_n \rangle} ,  \\
&(\text{\rm nI}_4) \, \langle x, y \mid x_2, \cdots, x_n \rangle =\langle x, y
\mid x_{j_2}, \cdots, x_{j_n} \rangle  \text{ for every permutation} \
(j_2, \cdots j_n) \\ & \qquad \qquad  \text{of} \ (2, \cdots, n) , \\
&(\text{\rm nI}_5) \,\ \langle x, x \mid x_2, x_3, \cdots, x_n \rangle = \langle x_2, x_2 \mid x, x_3, \cdots, x_n \rangle ,  \\
&(\text{\rm nI}_6) \,\ \langle \alpha x, y \mid x_2, \cdots, x_n\rangle = \alpha \langle x, y \mid
 x_2, \cdots, x_n \rangle , \\
&(\text{\rm nI}_7) \,\ \langle x+y, z \mid x_2, \cdots, x_n \rangle \le \langle x, z \mid x_2, \cdots,
x_n \rangle + \langle y, z \mid x_2, \cdots, x_n \rangle
\endalign$$
for all $\alpha \in \Bbb C$ and all $x, y, z, x_1, \cdots, x_n \in X$. The function
$\langle \cdot, \cdot \mid \cdot, \cdots, \cdot\rangle $
is called an {\it $n$-inner product on $X$}.
\enddefinition

The concept of $n$-inner product was introduced in \cite {1} for $n=2$ and in \cite{10} for $n\geq 2$.

For instance, on a given inner product space with inner product
$\langle ~,~ \rangle$, one can put an $n$-inner product by
defining $\langle x, y \mid x_2, \cdots, x_n \rangle$ to be
$$\operatorname{det} 
\pmatrix \langle x, y \rangle & \langle x, x_2 \rangle & \cdots &
\langle x, x_n \rangle \\ \langle x_2, y \rangle & \langle x_2,
x_2 \rangle & \cdots & \langle x_2, x_n \rangle \\ \vdots & \vdots
& \cdots & \vdots \\ \langle x_n, y \rangle & \langle x_n, x_2
\rangle & \cdots & \langle x_n, x_n \rangle  . \endpmatrix $$

\definition{Definition 3.2} An $n$-inner product and normed (respectively, Banach) space $X$ with $n$-inner product
$\langle \cdot, \cdot \mid \cdot, \cdots, \cdot\rangle_X $ and norm
$\| ~ \cdot ~ \|$ is called an {\it $n$-inner product normed (respectively, Banach) space}.
\enddefinition

For example, if $(X, \langle \cdot, \cdot \mid \cdot, \cdots, \cdot\rangle )$ is an $n$-inner product space, $\{a_1, \cdots, a_n\}$ is a fixed linearly independent set and $k>0$ then $\langle x, y \rangle: = k \displaystyle{\sum_{\{i_2, \cdots, i_n\}\subseteq\{1, \cdots, n\}}}\langle x, y \mid a_{i_2}, \cdots, a_{i_n} \rangle$ is an inner product on $X$ (see \cite{6}) and so $X$ equipped with the norm induced by this inner product is an $n$-inner product normed space.

In the rest of this section, assume that $X$ is an $n$-inner product Banach space with $n$-inner product $\langle \cdot, \cdot \mid \cdot, \cdots, \cdot\rangle_X $ and norm $||~ \cdot ~||$,
 and that $Y$ is an $n$-inner product Banach space with $n$-inner product
$\langle \cdot, \cdot \mid \cdot, \cdots, \cdot\rangle_Y$ and norm $\| ~ \cdot ~\|$.

\proclaim{Theorem 3.1}
 Let $f : X \rightarrow Y$ be a mapping with $f(0) = 0$
 for which there exists a function $\varphi
: X^{n+1} \rightarrow [0, \infty)$  such that
$$\align
\widetilde{\varphi}(x_0, \cdots, x_n) : =  \sum_{j=0}^{\infty}
\frac{1}{2^{j}}
\varphi(2^j x_0, & \cdots, 2^j x_n)   < \infty , \tag 3.i \\
 \|f(\mu x_0 + \mu x_1) - \mu f(x_0) - \mu f(x_1) \|
& \le \varphi(x_0, x_1, \underbrace{0, \cdots, 0}_{\text{$n-1$ times}}) ,  \tag 3.ii \\
| ~ \langle f(x_0), f(x_1) \mid f(x_2), \cdots, f(x_n)\rangle_Y &
- \langle x_0, x_1 \mid x_2, \cdots, x_n\rangle_X ~ | \\ &  \le
\varphi(x_0, \cdots, x_n)  \tag 3.iii
\endalign$$ for all $\mu \in \Bbb T^1$ and all $x_0,  \cdots, x_n \in X$.
 Then there exists a unique $\Bbb C$-linear $n$-inner product preserving mapping $U : X \rightarrow Y$ such that
$$\|f(x)- U(x)\| \le \frac{1}{2} \widetilde{\varphi}(x, x, \underbrace{0, \cdots, 0}_{\text{$n-1$ times}})
\tag 3.iv $$
for all $x \in X$.
\endproclaim

\demo{Proof}
By the same reasoning as in the proof of \cite{20, Theorem 2.1}, it follows from {\rm (3.i)} and {\rm (3.ii)} that
there exists a unique $\Bbb C$-linear mapping $U : X \rightarrow
Y$ satisfying {\rm (3.iv)}. The $\Bbb C$-linear mapping $U : X \rightarrow Y$ is given by
$$U(x) = \lim_{l \to \infty} \frac{1}{2^l} f(2^l x) \tag 3.1
$$ for all $x \in X$.

It follows from {\rm (3.iii)} that $$\align  | ~ \langle
\frac{1}{2^l} f(2^l x_0), \frac{1}{2^l} f(2^l x_1) \mid & \frac{1}{2^l} f(2^l x_2), \cdots,  \frac{1}{2^l} f(2^l x_n)\rangle_Y
 \\ & - \langle x_0, x_1 \mid x_2, \cdots, x_n\rangle_X ~ | \\ & =
 \frac{1}{2^{2nl}} | ~ \langle f(2^l x_0), f(2^l x_1) \mid f(2^l x_2),  \cdots, f(2^l x_n)\rangle_Y  \\ & \qquad  -
 \langle 2^l x_0, 2^l x_1 \mid 2^l x_2, \cdots,
2^l x_n\rangle_X ~ | \\ & \le \frac{1}{2^{2nl}}
 \varphi(2^l x_0, \cdots, 2^l x_n)  \le \frac{1}{2^{l}}
 \varphi(2^l x_0, \cdots, 2^l x_n) ,
\endalign $$
which tends to zero as $l \to \infty$ for all $x_0, \cdots, x_n \in X$ by {\rm (3.i)}. By {\rm (3.1)},
$$\align  \langle U(x_0), U(x_1) & \mid U(x_2), \cdots, U(x_n)\rangle_Y  \\ & = \lim_{l \to \infty} \langle \frac{1}{2^l}
f(2^l x_0), \frac{1}{2^l} f(2^l x_1) \mid \frac{1}{2^l} f(2^l
x_2), \cdots, \frac{1}{2^l} f(2^l x_n)\rangle_Y  \\ & = \langle
x_0, x_1 \mid x_2, \cdots, x_n\rangle_X \endalign $$ for all $x_0,
\cdots, x_n \in X$, as desired. \qed
\enddemo

\proclaim{Corollary 3.2}  Let $f : X \rightarrow Y$ be a
mapping with $f(0) = 0$ for which there exist constants
$\theta \ge 0$ and $p \in [0, 1)$ such that
$$\align \| f(\mu x_0 + \mu x_1)- \mu f(x_0) - \mu f(x_1)\| \le
\theta (||x_0||^p & + ||x_1||^p) , \\
| ~ \langle f(x_0), f(x_1) \mid f(x_2), \cdots, f(x_n)\rangle_Y -
\langle x_0, x_1 \mid x_2, \cdots, x_n\rangle_X ~ |  \le & \theta
\sum_{j=0}^n ||x_j||^p
\endalign$$ for all $\mu \in \Bbb T^1$ and all  $x_0, \cdots, x_n \in X$.
  Then there exists a unique $\Bbb C$-linear $n$-inner product preserving mapping  $U :  X
\rightarrow Y$ such that $$\|f(x)- U(x)\| \le \frac{2\theta}{2-2^p} ||x||^p $$
for all $x \in X$.
\endproclaim

\demo{Proof} Define $\varphi(x_0, \cdots, x_n) = \theta \sum_{j=0}^n ||x_j||^p $
to be Th.M. Rassias upper bound in the Cauchy--Rassias inequality,
and apply Theorem 3.1. \qed
\enddemo

\proclaim{Theorem 3.3}
 Let $f : X \rightarrow Y$ be a
mapping with $f(0)= 0$ for which there exists a function $\varphi
: (X \setminus \{0\})^{n+1} \rightarrow [0, \infty)$  such that
$$\align
\widetilde{\varphi}(x_0, \cdots, x_n) : = \sum_{j=0}^{\infty}
\frac{1}{3^j}
\varphi(3^j x_0, & \cdots,  3^j x_n)  < \infty , \tag 3.v \\
 \|2f(\frac{\mu x_0 + \mu x_1}{2}) -  \mu f(x_0) - \mu & f(x_1)\|
 \le \varphi(x_0, x_1, \underbrace{0, \cdots, 0}_{\text{$n-1$ times}}) ,  \tag 3.vi \\
| ~ \langle f(x_0), f(x_1) \mid f(x_2), \cdots, f(x_n)\rangle_Y &
- \langle x_0, x_1 \mid x_2, \cdots, x_n\rangle_X ~ | \\ & \le
\varphi(x_0, \cdots, x_n)  \tag 3.vii
\endalign$$ for all $\mu \in \Bbb T^1$ and all $x_0, \cdots, x_n \in X\setminus \{ 0 \}$.  Then there exists a unique
 $\Bbb C$-linear $n$-inner product preserving mapping $U :
X \rightarrow Y$ such that
$$\|f(x)- U(x)\| \le  \frac{1}{3}\big(\widetilde{\varphi}(x, -x, \underbrace{0, \cdots, 0}_{\text{$n-1$ times}}) +
\widetilde{\varphi}(-x, 3x, \underbrace{0, \cdots, 0}_{\text{$n-1$ times}})\big) \tag 3.viii $$
for all $x \in X \setminus \{0\}$.
\endproclaim

\demo{Proof} By the same reasoning as in the proof of \cite{20, Theorem 2.5},  it follows from {\rm (3.v)} and {\rm (3.vi)} that
there exists a unique $\Bbb C$-linear mapping $U : X \rightarrow
Y$ satisfying {\rm (3.viii)}. The $\Bbb C$-linear mapping $U : X
\rightarrow Y$ is given  by
$$U(x) = \lim_{l\to\infty}\frac{1}{3^l} f(3^l x) \tag 3.2 $$
for all $x \in X$.

It follows from {\rm (3.vii)} that $$\align  | ~ \langle
\frac{1}{3^l} f(3^l x_0), \frac{1}{3^l} f(3^l x_1) \mid & \frac{1}{3^l} f(3^l x_2), \cdots,
 \frac{1}{3^l} f(3^l x_n)\rangle_Y
 \\ & - \langle x_0, x_1 \mid x_2, \cdots, x_n\rangle_X ~ | \\ & =
 \frac{1}{3^{2nl}} | ~ \langle f(3^l x_0), f(3^l x_1) \mid f(3^l x_2),  \cdots, f(3^l x_n)\rangle_Y \\ & \qquad \quad -
 \langle 3^l x_0, 3^l x_1 \mid 3^l x_2, \cdots,
3^l x_n\rangle_X ~ | \\ & \le \frac{1}{3^{2nl}}
 \varphi(3^l x_0, \cdots, 3^l x_n) \le \frac{1}{3^{l}}
 \varphi(3^l x_0, \cdots, 3^l x_n) ,
\endalign $$
which tends to zero as $l \to \infty$ for all $x_0, \cdots, x_n \in X$ by {\rm (3.v)}. By {\rm (3.2)},
$$ \align \langle U(x_0), U(x_1) & \mid U(x_2), \cdots, U(x_n)\rangle_Y \\ & = \lim_{l \to \infty} \langle \frac{1}{3^l}
f(3^l x_0), \frac{1}{3^l} f(3^l x_1) \mid \frac{1}{3^l} f(3^l
x_2), \cdots, \frac{1}{3^l} f(3^l x_n)\rangle_Y \\ & = \langle
x_0, x_1 \mid x_2, \cdots, x_n\rangle_X \endalign  $$ for all
$x_0, \cdots, x_n \in X$, as desired. \qed
\enddemo

\proclaim{Corollary 3.4}  Let $f : X \rightarrow Y$ be a
mapping with $f(0) = 0$ for which there exist constants
$\theta \ge 0$ and $p \in [0, 1)$ such that
$$\align \|2f(\frac{\mu x_0 + \mu x_1}{2}) -  \mu f(x_0) - \mu f(x_1)\|  \le  \theta
(||x_0||^p & + ||x_1||^p) , \\
| ~ \langle f(x_0), f(x_1) \mid f(x_2), \cdots, f(x_n)\rangle_Y -
\langle x_0, x_1 \mid x_2, \cdots, x_n\rangle_X ~ | \le & \theta
\sum_{j=0}^n ||x_j||^p
\endalign$$ for all $\mu \in \Bbb T^1$ and all  $x_0, \cdots, x_n \in X\setminus \{ 0 \}$.
  Then there exists a unique $\Bbb C$-linear $n$-inner product preserving mapping $U :  X
\rightarrow Y$ such that $$\|f(x)- U(x)\| \le \frac{3+3^p}{3-3^p} \theta ||x||^p $$
for all $x \in X$.
\endproclaim

\demo{Proof} Define $\varphi(x_0, \cdots, x_n) = \theta \sum_{j=0}^n ||x_j||^p$,
and apply Theorem 3.3. \qed
\enddemo

\proclaim{Theorem 3.5}
 Let $f : X \rightarrow Y$ be a
mapping with $f(0)= 0$ for which there exists a function $\varphi
: (X \setminus \{0\})^{n+1} \rightarrow [0, \infty)$  such that
$$\align
 \sum_{j=0}^{\infty} 3^{2nj}
\varphi(\frac{x_0}{3^j}, \cdots, \frac{x_n}{3^j}) &  < \infty , \tag 3.ix \\
 \|2f(\frac{\mu x_0 +\mu x_1}{2}) - \mu  f(x_0) - \mu f(x_1)\|
& \le \varphi(x_0, x_1, \underbrace{0, \cdots, 0}_{\text{$n-1$ times}}) ,  \tag 3.x \\
| ~ \langle f(x_0), f(x_1) \mid f(x_2), \cdots, f(x_n)\rangle_Y &
- \langle x_0, x_1 \mid x_2, \cdots, x_n\rangle_X ~ |  \\ & \le
\varphi(x_0, \cdots, x_n)   \tag 3.xi
\endalign$$ for all $\mu \in \Bbb T^1$ and all  $x_0, \cdots, x_n \in X \setminus \{ 0 \}$.
 Then there exists a unique $\Bbb C$-linear $n$-inner product preserving mapping $U :
X \rightarrow Y$ such that
$$\|f(x)- U(x)\| \le  \widetilde{\varphi}(\frac{x}{3}, - \frac{x}{3}, \underbrace{0, \cdots, 0}_{\text{$n-1$ times}}) +
\widetilde{\varphi}(-\frac{x}{3}, x, \underbrace{0, \cdots, 0}_{\text{$n-1$ times}}) \tag 3.xii $$
for all $x \in X \setminus \{0\}$, where $$ \widetilde{\varphi}(x_0, \cdots, x_n) : =  \sum_{j=0}^{\infty} 3^{j}
\varphi(\frac{x_0}{3^j}, \cdots, \frac{x_n}{3^j}) $$ for all $x_0, \cdots, x_n \in X$.
\endproclaim

\demo{Proof} Note that $$3^{j} \varphi(\frac{x_0}{3^j}, \cdots, \frac{x_n}{3^j}) \le
3^{2nj} \varphi(\frac{x_0}{3^j}, \cdots, \frac{x_n}{3^j}) \tag 3.3$$
for all $x_0, \cdots, x_n \in X$ and all positive integers $j$.
By the same reasoning as in the proof of Theorem 2.5,
 it follows from {\rm (3.ix)}, {\rm (3.3)} and {\rm (3.x)} that
there exists a unique $\Bbb C$-linear mapping $U : X \rightarrow
Y$ satisfying {\rm (3.xii)}. The $\Bbb C$-linear mapping $U : X
\rightarrow Y$ is given  by
$$U(x) = \lim_{l\to\infty}3^l f(\frac{x}{3^l}) \tag 3.4 $$
for all $x \in X$.

It follows from  {\rm (3.xi)} that $$\align  | ~ \langle 3^l
f(\frac{x_0}{3^l}), 3^l f(\frac{x_1}{3^l}) \mid & ~ 3^l
f(\frac{x_2}{3^l}), \cdots, 3^l f(\frac{x_n}{3^l})\rangle_Y
 \\ & - \langle x_0, x_1 \mid x_2, \cdots, x_n\rangle_X ~ | \\ & =
 3^{2nl} | ~ \langle f(\frac{x_0}{3^l }), f(\frac{x_1}{3^l }) \mid f(\frac{x_2}{3^l}),  \cdots,
f(\frac{x_n}{3^l})\rangle_Y  \\ & \qquad  \quad -
 \langle \frac{x_0}{3^l}, \frac{x_1}{3^l} \mid \frac{x_2}{3^l}, \cdots,
\frac{x_n}{3^l} \rangle_X ~ | \\ & \le 3^{2nl}
 \varphi(\frac{x_0}{3^l}, \cdots, \frac{x_n}{3^l})  ,
\endalign $$
which tends to zero as $l \to \infty$ for all $x_0, \cdots, x_n \in X$ by {\rm (3.ix)}. By {\rm (3.4)},
$$\align  \langle U(x_0), U(x_1) & \mid U(x_2), \cdots, U(x_n)\rangle_Y \\ & = \lim_{l \to \infty} \langle 3^l
f(\frac{x_0}{3^l}), 3^l f(\frac{x_1}{3^l}) \mid 3^l
f(\frac{x_2}{3^l}), \cdots, 3^l f(\frac{x_n}{3^l})\rangle_Y \\ & =
\langle x_0, x_1 \mid x_2, \cdots, x_n\rangle_X \endalign $$ for
all $x_0, \cdots, x_n \in X$, as desired. \qed
\enddemo

\proclaim{Corollary 3.6}  Let $f : X \rightarrow Y$ be a
mapping with $f(0) = 0$ for which there exist constants
$\theta \ge 0$ and $p \in (2n, \infty)$ such that
$$\align \|2f(\frac{\mu x_0 + \mu x_1}{2}) -  \mu f(x_0) - \mu f(x_1)\| \le  \theta
(||x_0||^p & + ||x_1||^p) , \\
| ~ \langle f(x_0), f(x_1) \mid f(x_2), \cdots, f(x_n)\rangle_Y -
\langle x_0, x_1 \mid x_2, \cdots, x_n\rangle_X ~ | \le & \theta
\sum_{j=0}^n ||x_j||^p
\endalign$$ for all $\mu \in \Bbb T^1$ and all  $x_0, \cdots, x_n \in X\setminus \{ 0 \}$.
  Then there exists a unique $\Bbb C$-linear $n$-inner product preserving mapping $U :  X
\rightarrow Y$ such that $$\|f(x)- U(x)\| \le \frac{3^p+3}{3^p-3} \theta ||x||^p $$
for all $x \in X\setminus \{ 0 \}$.
\endproclaim

\demo{Proof} Define $\varphi(x_0, \cdots, x_n) = \theta \sum_{j=0}^n ||x_j||^p$,
and apply Theorem 3.5. \qed
\enddemo

{\bf Acknowledgment.} We sincerely thank Professor Themistocles M. Rassias for his valuable suggestions.

\enddocument